\documentclass[11pt,a4paper]{article}
\usepackage{amsmath,amsthm,amsfonts, graphicx, amssymb,iwona}
\usepackage[numbers,sort&compress]{natbib}
\usepackage[lmargin=35mm,rmargin=35mm,tmargin=30mm,bmargin=30mm]{geometry}
\renewcommand{\baselinestretch}{1.25}
\usepackage[colorlinks=true,citecolor=black,linkcolor=black,urlcolor=blue]{hyperref}
\newcommand{\msn}[1]{MR:\,\href{http://www.ams.org/mathscinet-getitem?mr=MR#1}{#1}}
\newcommand{\doi}[1]{\href{http://dx.doi.org/#1}{\texttt{doi:#1}}}

\newcommand{\urlprefix}{}
\theoremstyle{plain}
\newtheorem{theorem}{Theorem}
\newtheorem{lemma}{Lemma}
\allowdisplaybreaks
\renewcommand{\geq}{\geqslant}
\renewcommand{\leq}{\leqslant}

\newcommand{\floor}[1]{\ensuremath{\protect\lfloor#1\rfloor}}

\newcommand{\half}{\ensuremath{\protect\tfrac{1}{2}}}

\begin{document}

\title{\vspace*{-3ex}\bf\LARGE On Multiplicative Sidon Sets}

\author{David Wakeham\thanks{Department of Mathematics and
  Statistics, The University of Melbourne, Melbourne,
  Australia (\texttt{d.wakeham@student.unimelb.edu.au}).} \qquad
David~R.~Wood\thanks{Department of Mathematics and Statistics, The
  University of Melbourne, Melbourne, Australia
  (\texttt{woodd@unimelb.edu.au}). Supported by a QEII Research
  Fellowship from the Australian Research Council.} }

\date{\today}

\maketitle

\vspace*{-2ex}
\begin{abstract}
  \noindent Fix integers $b>a\geq1$ with $g:=\gcd(a,b)$. A set
  $S\subseteq\mathbb{N}$ is \emph{$\{a,b\}$-multiplicative} if $ax\neq
  by$ for all $x,y\in S$.  For all $n$, we determine an
  $\{a,b\}$-multiplicative set with maximum cardinality in $[n]$, and
  conclude that the maximum density of an $\{a,b\}$-multiplicative set
  is  $\frac{b}{b+g}$.
  For $A, B \subseteq \mathbb{N}$, a set $S\subseteq\mathbb{N}$ is
  \emph{$\{A,B\}$-multiplicative} if $ax=by$ implies $a = b$ and $x =
  y$ for all $a\in A$ and $b\in B$, and $x,y\in S$.
  For $1 < a < b < c$ and $a, b, c$ coprime, we give an $O(1)$ time
  algorithm to approximate the maximum density of an $\{\{a\},\{b,
  c\}\}$-multiplicative set to arbitrary given precision.
  %We generalise this notion to $\{A,
%  B\}$-multiplicativity, $A, B \subseteq \mathbb{N}$ and approximate
%  the maximum density in the simplest nontrivial case where $|A| = 1$,
%  $|B| = 2$.
\end{abstract}

\bigskip
\section{Introduction}
\citet{Erdos69, Erdos38, Erdos6869} defined a set $S\subseteq
\mathbb{N}$ to be \emph{multiplicative Sidon}\footnote{Additive Sidon
  sets have been more widely studied; see the classical papers
  \citep{Sidon32, Singer38, ET41} and the survey by
  \citet{OBryant-EJC04}. Let $\mathbb{N}:=\{1,2,\dots\}$,
  $\mathbb{N}_0 = \mathbb{N} \cup \{0\}$ and
  $[n]:=\{1,2,\dots,n\}$.} if $ab=cd$ implies $\{a,b\}=\{c,d\}$ for
all $a,b,c,d\in S$; see \citep{Sarkozy01,Ruzsa-AMH06,Ruzsa-JNT99}. In a
similar direction, \citet{Wang} defined a set $S\subseteq\mathbb{N}$
to be \emph{double-free} if $x\neq 2y$ for all $x,y\in S$, and proved
that the maximum density of a double-free set is $\frac{2}{3}$; see
\citep{AABBJPS-DM95} for related results. Here the \emph{density} of
$S\subseteq\mathbb{N}$ is $$\lim_{n\rightarrow\infty}\frac{|S\cap
  [n]|}{n}\enspace.$$ Motivated by some questions in graph colouring,
\citet{PorWood-Comb09} generalised the notion of double-free sets as
follows. For $k\in\mathbb{N}$, a set $S\subseteq\mathbb{N}$ is
\emph{$k$-multiplicative} (\emph{Sidon}) if $ax=by$ implies $a=b$ and
$x=y$ for all $a,b\in[k]$ and $x,y\in S$. \citet{PorWood-Comb09}
proved that the maximum density of a $k$-multiplicative set is
$\Theta(\frac{1}{\log k})$.

Here we study the following alternative generalisation of double-free
sets. 
For distinct $a,b\in\mathbb{N}$, a set $S\subseteq\mathbb{N}$ is
\emph{$\{a,b\}$-multiplicative} if $ax\neq by$ for all $x,y\in S$.
Our first result is to determine the maximum density of an
$\{a,b\}$-multiplicative set. Assume that $a<b$ throughout.

Say $x\in\mathbb{N}$ is an \emph{$i$-th subpower} of $b$ if $x=b^iy$
for some $y\not\equiv 0\pmod{b}$.
If $x$ is an $i$-th subpower of $b$ for some even/odd $i$ then $x$ is
an \emph{even/odd} subpower of $b$.
We prove the following result:

\begin{theorem}
  Fix integers $b>a\geq1$. Let $g:=\gcd(a,b)$. Then for every integer
  $n\in\mathbb{N}$, the even subpowers of $\frac{b}{g}$ in $[n]$ are an
  $\{a,b\}$-multiplicative set in $[n]$ with maximum cardinality. And the even
  subpowers of $\frac{b}{g}$ are an $\{a,b\}$-multiplicative set with
  density $\frac{b}{b+g}$, which is maximum.
\end{theorem}

\noindent Note that if $g=a$ then $b\geq 2g$ and $b+g\leq\frac{3}{2}b$, and
 if $g<a$ then $a\geq 2g$ and $b+g\leq b+a<\frac{3}{2}b$. In both
 cases the density bound $\frac{b}{b+g}$ in Theorem~1 is at least
 $\frac{2}{3}$, which is the bound obtained by \citet{Wang} for the
 $a=1$ and $b=2$ case.

We propose a further generalisation of double-free sets.
Let $A,B\subseteq\mathbb{N}$.
Say $S\subseteq\mathbb{N}$ is \emph{$\{A,B\}$-multiplicative} if
$ax=by$ implies $a = b$ and $x = y$ for all $a\in A$ and $b\in B$, and
$x,y\in S$.
One case is easily dealt with. 
%For some prime number $b$, let $A:=[b-1]$ and $B:=\{b\}$.
%Then $\gcd(a,b)=1$ for all $a\in A$.
%Thus the even subpowers of $b$ are $\{A,B\}$-multiplicative, and have
%maximum density.
If $B := \{b\}$ and $b$ is coprime to each element of $A$, and there
is some element $a \in A$ such that $a < b$, then, by the reasoning
above, the even subpowers of $b$ form an $\{A, B\}$-multiplicative set
of (maximum) density $\frac{b}{b + 1}$.

The simplest nontrivial case (not covered by Theorem 1) is $\{A,
B\}$-multiplicativity for $A = \{a\}$, $B = \{b, c\}$, $1 < a < b <
c$, with $a, b, c$ pairwise coprime.
We have the following theorem:

\begin{theorem}
Consider $a, b, c \in \mathbb{N}$ pairwise coprime, $1 < a < b < c$.
For all fixed $\epsilon > 0$, there is an $\mathcal{O}(1)$ time
algorithm that computes the maximum density of an $\{\{a\}, \{b,
c\}\}$-multiplicative set to within $\epsilon$.
\end{theorem}

\section{Proof of Theorem 1}

First suppose that $\gcd(a,b)=1$.  Let $T$ be the set of even
subpowers of $b$. We now prove that $T$ is an $\{a,b\}$-multiplicative
set with maximum density. In fact, for all $[n]$, we prove that
$T_n:=T\cap[n]$ has maximum cardinality out of all
$\{a,b\}$-multiplicative sets contained in $[n]$.

The key to our proof is to model the problem using a directed graph.
Let $G$ be the directed graph with $V(G):=[n]$ where $xy\in E(G)$
whenever $bx=ay$ (implying $x<y$). Thus $S\subseteq[n]$ is
$\{a,b\}$-multiplicative if and only if $S$ is an independent set in
$G$. If $xyz$ is a directed path in $G$, then $x=\frac{a}{b}y$ and
$z=\frac{b}{a}y$. Thus each vertex $y$ has indegree and outdegree at
most $1$. Since $xy\in E(G)$ implies $x<y$, $G$ contains no directed
cycles. Thus $G$ is a collection of disjoint directed paths.  Hence a
maximum independent set in $G$ is obtained by taking all the vertices
at even distance from the source vertices\footnote{Note that this is not
necessarily the only maximum independent set---for a path component
with odd length, we may take the vertices at odd distance from the
source of this path. This observation readily leads to a
characterisation of all maximum independent sets in $G$, and thus of
all $\{a,b\}$-multiplicative sets in $[n]$ with maximum
cardinality. Details omitted.}, where a vertex $y$ is a source (indegree
0) if and only if $\frac{a}{b}y$ is not an integer; that is, if
$y\not\equiv0\pmod{b}$.

We now prove that the vertices at distance $d$ from a source vertex
are precisely the $d$-th subpowers of $b$. We proceed by induction on
$d\geq 0$.  Each vertex $y$ of $G$ has an incoming edge if and only if
$\frac{a}{b}y\in\mathbb{N}$, which occurs if and only if
$y\equiv0\pmod{b}$ since $\gcd(a,b)=1$. Thus the source vertices of
$G$ are precisely the $0$-th subpowers of $b$. This proves the $d=0$
case of the induction hypothesis. Now consider a vertex $y$ at
distance $d$ from a source vertex. Thus $y=\frac{b}{a}x$ for some
vertex $x$ at distance $d-1$ from a source vertex. By induction, $x$
is a $(d-1)$-th subpower of $b$. That is, $x=b^{d-1}z$ for some
$z\not\equiv0\pmod{b}$. Thus $y=b^d\frac{z}{a}$, which, since
$\gcd(a,b)=1$,  implies that
$\frac{z}{a}$ is an integer. Hence $\frac{z}{a}\not\equiv0\pmod{b}$ and $y$ is a
$d$-th subpower of $b$, as claimed.

This proves that the even subpowers of $b$ form a maximum independent
set in $G$. That is, $T_n$ is an $\{a,b\}$-multiplicative set of
maximum cardinality in $[n]$. To illustrate this proof, the following
table shows two examples of the graph $G$ with $b=3$.  Observe that
the $i$-th row consists of the $i$-th subpowers of $3$ regardless of
$a$.

{\medskip\footnotesize \setlength{\tabcolsep}{4.2pt}
\begin{center}
 \begin{tabular}{|ccccccccc|cccccccccccc|}
    \hline
    \multicolumn{9}{|c|}{$a=1$ and $b=3$}&
    \multicolumn{12}{c|}{$a=2$ and $b=3$}\\
    \hline
    1&2&4&5&7&8&10&11&$\cdots$&
    1&2&4&5&7&8&10&11&13&14&16&$\cdots$\\
%%%%%%%%%%%%
    $\downarrow$&$\downarrow$&$\downarrow$&$\downarrow$&$\downarrow$&$\downarrow$&$\downarrow$&$\downarrow$&~&
    &$\downarrow$&$\downarrow$&&&$\downarrow$&$\downarrow$&&&$\downarrow$&$\downarrow$&~\\
%%%%%%%%%%%%
    3&6&12&15&21&24&30&33 &$\cdots$&
    & 3&6&&&12&15&&&21&24 &$\cdots$\\
%%%%%%%%%%%%
    $\downarrow$&$\downarrow$&$\downarrow$&$\downarrow$&$\downarrow$&$\downarrow$&$\downarrow$&$\downarrow$&~&
    &&$\downarrow$&&&$\downarrow$&&&&&$\downarrow$&~\\
%%%%%%%%%%%%
    9&18&36&45&63&72&90&99& $\cdots$&
    && 9&&&18&&&&&36& $\cdots$\\
%%%%%%%%%%%%
    $\downarrow$&$\downarrow$&$\downarrow$&$\downarrow$&$\downarrow$&$\downarrow$&$\downarrow$&$\downarrow$&~&
    &&&&&$\downarrow$&&&&&$\downarrow$&~\\
%%%%%%%%%%%%
    27&48&108&135&189&216&270&297&$\cdots$
    &&&&&&27&&&&&48&$\cdots$\\
    $\vdots$&$\vdots$&$\vdots$&$\vdots$&$\vdots$&$\vdots$&$\vdots$&$\vdots$&
    &&&&&&&&&&&$\vdots$&\\
    \hline
  \end{tabular}\end{center}\medskip}

%CAN WE BE MORE PRECISE ABOUT $|T_n|$??? See Wang's paper.

We now bound $|T_n|$ from above. Observe that
$$T_n=\left\{b^{2i}y: 0\leq i\leq\half\log_b n, \, 1\leq y\leq \frac{n}{b^{2i}},\,
  y\not\equiv0 \!\!\!\!\pmod{b}\right\}\enspace.$$ Thus
\begin{align*}
  |T_n| & \leq \sum_{i=0}^{\floor{(\log_b n)/2}}
  \left\lceil\frac{b-1}{b}\frac{n}{b^{2i}}\right\rceil\\
%%%%%%%%%%%%%
  & \leq 1+ \half(\log_b  n)+\frac{(b-1)n}{b}\sum_{i\geq0}\frac{1}{b^{2i}}\\
%%%%%%%%%%%%%
  & \leq 1+ \half(\log_b  n)+\frac{(b-1)n}{b}\,\frac{b^2}{b^2-1}\\
%%%%%%%%%%%%%
  & = 1+\half(\log_b n)+\frac{bn}{b+1}\enspace.
\end{align*}

We now bound $|T_n|$ from below. Observe that
$$T_n=[n]\setminus\left\{b^{2i+1}y: 0\leq i\leq\half((\log_b n)-1), \, 1\leq y\leq \frac{n}{b^{2i+1}},\,
  y\not\equiv0\pmod{b}\right\}\enspace.$$ Thus
\begin{align*}
  |T_n| & \geq n-\sum_{i=0}^{\floor{((\log_b n)-1)/2}}
  \left\lceil\frac{b-1}{b}\frac{n}{b^{2i+1}}\right\rceil\\
  & \geq n-\half((\log_b n)+1)-\frac{(b-1)n}{b^2}\sum_{i\geq  0}\frac{1}{b^{2i}}\\
  & \geq n-\half((\log_b n)+1)-\frac{(b-1)n}{b^2}\,\frac{b^2}{b^2-1}\\
  & = n-\half((\log_b n)+1)-\frac{n}{b+1}\\
  & = \frac{bn}{b+1}-\half((\log_b n)+1)\enspace.
\end{align*}

These upper and lower bounds on $|T_n|$ imply that
$$|T_n|= \frac{bn}{b+1}+\Theta(\log_b n)\enspace.$$
Hence the density of $T$
is $\frac{b}{b+1}$, and because $T_n$ is optimal for each $n$, no
$\{a,b\}$-multiplicative set has density greater than $\frac{b}{b+1}$.

We now drop the assumption that $\gcd(a,b)=1$. Let $g:=\gcd(a,b)$.
Since $ax=by$ if and only if $\frac{a}{g}x=\frac{b}{g}y$, a set $S$ is
$\{a,b\}$-multiplicative if and only if $S$ is
$\{\frac{a}{g},\frac{b}{g}\}$-multiplicative. Since
$\frac{b/g}{b/g+1}=\frac{b}{b+g}$, the theorem is proved.

\section{Proof of Theorem 2}

Fix $A = \{a\}$ and $B = \{b, c\}$, where $1 < a < b < c$, and $a, b, c$
pairwise coprime.
For convenience, we use the infinite graph $G$ with vertex set
$\mathbb{N}$ and edge set
\[
E(G) = \{xy : bx = ay \text{ or } cx = ay, x, y \in \mathbb{N}\}.
\]
Let $G_n$ denote the subgraph of $G$ induced by the vertex set $[n]$.
Let $\delta$ be the maximum density of an $\{\{a\}, \{b,
c\}\}$-multiplicative set.
Then
\[
\delta = \lim_{n \to \infty} \frac{\alpha(G_n)}{n} \enspace ,
\]
where $\alpha(G_n)$ is the size of a maximum independent set in
$G_n$.

The infinite graph $G$ has components $C_{p, q}$ with
vertex set $$V(C_{p, q}) = \{a^{p - x - y}b^xc^yq : x, y \in
\mathbb{N}_0\}$$ for all $p \in \mathbb{N}_0$, $q \in \mathbb{N}$, and
$q$ not divisible by $a$, $b$, or $c$.
Note that each $C_{p, q}$ is finite.
Define $p$ as the \emph{height} of the component, and subsets of
constant $x + y$ as \emph{rows}.
Note that the maximum and minimum vertices in $C_{p,q}$ are $c^pq$ and
$a^pq$ respectively.
The first few components of $G$ for $a = 2$, $b = 3$, and $c = 5$ are shown
below:
\vspace{10pt}
\begin{center}
  \includegraphics[scale = 0.35]{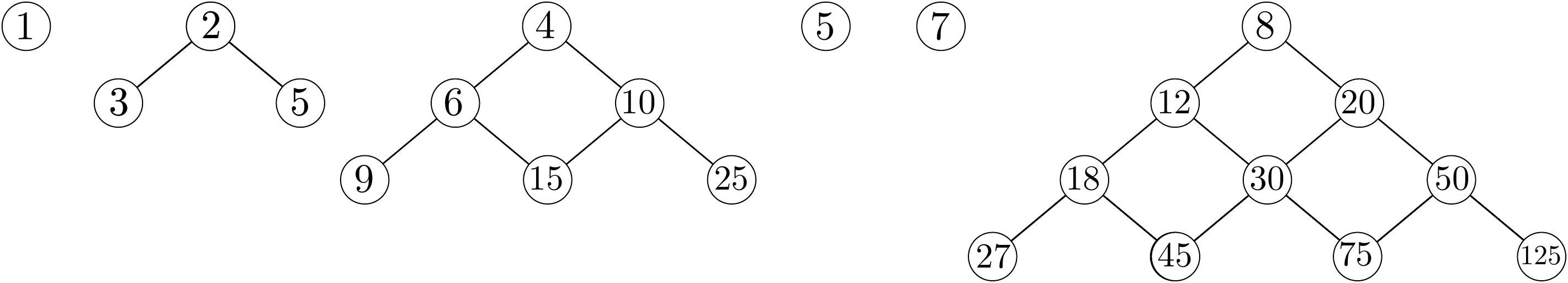}
\end{center}

For $a, b, c$ as above and fixed $\epsilon > 0$, let $d$ be a
non-negative integer $d \in \mathbb{N}_0$, to be specified
later.
Basically, $d$ is a cutoff height which allows us to partition the
components of $G_n$ into three types, for any given $n \in
\mathbb{N}$.
The first are \emph{complete} components $C_{p, q}$ where $n > c^pq$.
The second are \emph{small} incomplete components $S_{p, q}$ where $p
\leq d$ and $a^pq \leq n < c^pq$.
The third are \emph{large} incomplete components $L_{p, q}$ with $p > d$
and $a^pq \leq n < c^pq$.

Let $\alpha_T(G_n)$ denote the size of a maximum independent set in
the components of type $T$ in $G_n$, for $T \in \{C, S, L\}$.
We clearly have
$$\alpha(G_n) = \alpha_C(G_n) + \alpha_S(G_n) + \alpha_L(G_n).$$
Thus,
$$\delta = \lim_{n\to\infty} \frac{\alpha_C(G_n)}{n} +
\lim_{n\to\infty} \frac{\alpha_S(G_n)}{n} + \lim_{n\to\infty}
\frac{\alpha_L(G_n)}{n} =
\delta_C + \delta_S + \delta_L$$
where
\[
\delta_T = \lim_{n\to\infty} \frac{\alpha_T(G_n)}{n}.
\]
We determine $\delta_C$ and $\delta_S$ below, and show that, for any
$\epsilon > 0$, we can choose $d$ so that $\delta_L < \epsilon$.
Hence, we can calculate $\delta$ to arbitrary precision.

\subsection*{Complete components}

We require the following lemma about independent sets in
grid-like graphs by \citet{caszim}.

\begin{lemma}
Define a graph $H$ by $V(H) := \mathbb{N}_0 \times \mathbb{N}_0$ and
\[
E(H) := \{\{\mathbf{v}, \mathbf{w}\} : \mathbf{v}, \mathbf{w} \in
V(H), |v_1 - w_1| + |v_2 - w_2| = 1\}.
% = \{\{(m_1, n_1), (m_2, n_2)\} : |m_1 - m_2| + |n_1 - n_2| = 1\}.
\]
Suppose that $F$ is a finite subgraph of $H$ such that $(x, y) \in
V(F)$ implies $(x - 1, y) \in V(F)$ unless $x = 0$, and
$(x, y - 1) \in V(F)$ unless $y = 0$.
Then one of the sets
\begin{align*}
O & := \{(x, y) \in V(F) : x + y \text{\emph{ is odd}}\} \text{ or}\\
E &:= \{(x, y) \in V(F) : x + y \text{\emph{ is even}}\}
\end{align*}
is a maximum independent set in $F$.
\end{lemma}

Now, consider a complete component $C_{p, 1}$ of $G_n$.
Note that every complete component $C_{p, q}$ of height $q$ is
isomorphic to $C_{p, 1}$, and can be obtained by multiplying each
vertex by $q$.
Thus, we call $C_{p, q}$ a $q$-\emph{copy} of $C_{p, 1}$.
In general, we use this terminology for isomorphic components of any
type obtained by multiplying each vertex by $q$.

Observe that we can apply Lemma 1 to $C_{p, 1}$, since it is
isomorphic to a subgraph of $H$ with the required properties.
Define a function $\varphi: V(C_{p, 1}) \to \mathbb{N}_0 \times
\mathbb{N}_0$ by
\[
\varphi(a^{p-x-y}b^xc^y) = (x, y).
\]
If $a^{p-x-y}b^xc^y$ is adjacent to $a^{p-x'-y'}b^{x'}c^{y'}$, then
$|x - x'| + |y - y'| = 1$ since they must differ by a factor of $b/a$
or $c/a$.
Thus, since $\varphi$ is injective, it defines an isomorphism from
$C_{p, 1}$ to a subgraph of $H$.
Assume $a^{p-x-y}b^{x}c^{y} \in V(C_{p, 1})$.
Then $a^{p -x-y + 1}b^{x - 1}c^{y} \in V(C_{p, 1})$ unless $x = 0$,
and similarly $a^{p -x-y + 1}b^{x}c^{y - 1} \in V(C_{p, 1})$ unless $y
= 0$.
Under $\varphi$, these are clearly equivalent to the conditions
required for Lemma 1.

Hence, by Lemma 1 and the definition of $\varphi$, a maximum
independent set in $C_{p, 1}$ is given by choosing all rows with $x +
y$ even, or all rows with $x + y$ odd.
In fact, it is clear that a maximum independent set is obtained by
choosing the bottom row first, then alternating between remaining rows.
Thus, if $p = 2i - 1$, then $\alpha(C_{p, 1}) = i(i+1)$.
If $p = 2i$, then $\alpha(C_{p, 1}) = (i+1)^2$.
Since the largest vertex in such a component is $c^p$, we must have $p
\leq \log_cn$ for the component $C_{p, 1}$ to be complete.
Hence, the maximum height of a complete component is $\lfloor\log_cn\rfloor$.

Now we multiply by the number of components of height $p$ that are
complete.
For a given $p$, we require $1 \leq q \leq nc^{-p}$.
Since the density of numbers not divisible by $a$, $b$, or $c$ is
\[
\frac{(a-1)(b-1)(c-1)}{abc},
\]
the number of components of height $p$ in $G_n$ is
$$\frac{(a-1)(b-1)(c-1) n}{c^{p}abc} + o(n)\,.$$
Let $M(n) = \frac{1}{2}\lfloor \log_cn\rfloor$.
The total number of vertices in a maximum independent set in complete
components is therefore
$$\alpha_C(G_n) =
\frac{(a - 1)(b-1)(c-1)n}{abc}\sum_{i=0}^{M(n)}\left[\frac{i(i+1)}{c^{2i-1}}
  +\frac{(i+1)^2}{c^{2i}}\right] + o(n).$$
Thus, the density contribution is
\begin{align*}
\delta_C = \lim_{n\to\infty} \frac{\alpha_C(G_n)}{n} & =
\lim_{n\to\infty}\frac{1}{n}\cdot\frac{(a - 1)(b - 1)(c - 1)n}{abc}\sum_{i=0}^{M(n)}\left[\frac{i(i+1)}{c^{2i-1}}
  + \frac{(i+1)^2}{c^{2i}}\right] \\
& = \frac{(a - 1)(b - 1)(c - 1)}{abc}\sum_{i=0}^{\infty}\left[\frac{i(i+1)}{c^{2i-1}}
  + \frac{(i+1)^2}{c^{2i}}\right] \\
& = \frac{(a-1)(b-1)(c-1)}{abc} \cdot \frac{c^4}{(c-1)^3(c+1)} \\
& = \frac{(a-1)(b-1)c^3}{ab(c-1)^2(c+1)}.
\end{align*}
%For example, the density of complete components for $a = 2, b = 3, c =
%5$ is
%\[
%\frac{2\cdot 5^3}{2\cdot 3\cdot 4^2 \cdot 6}  = \frac{125}{288} \approx 0.434.
%\]

\subsection*{Small incomplete components}

Now we consider the small incomplete components.
Let $C_{p, 1}[r]$ be the subgraph of $C_{p, 1}$ induced by $[r]$.
Define
$$f(p, r) := \alpha(C_{p, 1}[r])$$
for $r \in \mathbb{N}$.
We can calculate all $f$ for $p \leq d$ in $\mathcal{O}(c^d)$ time
with a computer, again using Lemma 1.
(In fact, these components have bounded size, so any exponential time
maximum independent set algorithm runs in $O(1)$ time.)
Note that $C_{p,q}[n]$ is a $q$-copy of $C_{p, 1}[\lfloor
n/q\rfloor]$, and therefore $\alpha(C_{p, q}[n]) = f(p, \lfloor n/q
\rfloor)$.
So we can find the size of maximum independent sets in the small
components using the $f$'s.

More precisely, given $p \leq d$ and $n$, for how many values of $q$
is $C_{p, q}[n]$ a $q$-copy of $C_{p, 1}[r]$, where $r = \lfloor n/q\rfloor$?
First note that
$$\frac{n}{r+1} < q \leq \frac{n}{r}.$$
Thus, there are
$$\frac{(a-1)(b-1)(c-1)n}{abc}\left(\frac{1}{r} - \frac{1}{r+1}\right)
+ o(n)=
\frac{(a-1)(b-1)(c-1)n}{abcr(r+1)} + o(n)$$
$q$-copies of $C_{p, 1}[r]$.
%This expression is correct in the limit as $n \to \infty$.
The only restriction on $r$ is that $a^p \leq r \leq c^p - 1$.
Hence, the size of a maximum independent set in components of type $S$
is
$$\sum_{p = 0}^d \sum_{r = a^p}^{c^p - 1}
\frac{(a-1)(b-1)(c-1)n}{abcr(r+1)}f(p, r) + o(n).$$
As $n \to \infty$, the density contribution of small components is therefore
\begin{align*}
\delta_S & = \lim_{n \to \infty}\frac{1}{n} \sum_{p = 0}^d
\sum_{r = a^p}^{c^p - 1} \frac{(a-1)(b-1)(c-1)n}{abcr(r+1)}f(p, r) \\
& = \sum_{p = 0}^d \sum_{r = a^p}^{c^p - 1}
\frac{(a-1)(b-1)(c-1)}{abcr(r+1)}f(p, r).
\end{align*}
Since $a$, $b$, and $c$ are constants and $d$ is bounded by a function
of $a$, $b$, and $c$, $\delta_S$ can be computed in $\mathcal{O}(1)$
time. 

\subsection*{Large incomplete components}

Finally, we show that we can choose $d$ so that the density of a
maximum independent set in components of type $L$ is less than
$\epsilon$.
For large components, 
$$p > d \,\,\text{ and }\,\, a^pq \leq n < c^pq.$$
The latter implies $c^{-p}n < q \leq a^{-p}n$.
From the density of $q$, the number of large incomplete components
$L_{p, q}$ for a given $p > d$ is
$$\frac{(a-1)(b-1)(c-1)n}{abc}\left(\frac{1}{a^p} -
  \frac{1}{c^p}\right) + o(n).$$
Since there are less than $p^2$ vertices in a component of height $p$, 
\begin{align*}
\alpha_L(G_n) & \leq \sum_{p=d}^\infty
p^2\cdot\frac{(a-1)(b-1)(c-1)n}{abc}\left(\frac{1}{a^p} - \frac{1}{c^p}\right) \\
& \leq \frac{(a-1)(b-1)(c-1)n}{abc}\sum_{p=d}^\infty \frac{p^2}{a^p} \\
& = \frac{(a-1)(b-1)(c-1)n}{abc}\cdot \frac{a^{1-d}((a-1)^2d^2 + 2(a - 1)d
  + a + 1)}{(a-1)^3} \\
& \leq \frac{(b-1)(c-1)n}{bc}\cdot a^{-d/2}
\end{align*}
where the last inequality holds for $d \geq 22$.
Define $\beta := (b-1)(c-1)/bc$.
Hence,
$$\delta_L = \lim_{n\to\infty}\frac{\alpha_L(G_n)}{n} \leq
\lim_{n\to\infty} \frac{1}{n}\cdot \beta n \cdot a^{-d/2} = \beta a^{-d/2}.$$
So, to obtain a precision of $\epsilon$ in the approximation $\delta
\approx \delta_C + \delta_S$, we pick
$$d = \max\{2\log_{a}(\beta/\epsilon), 22\}$$
which is a function of $a, b, c$, and $\epsilon$.
This completes the proof of Theorem 2.

%In practice, this lower bound on $d$ is very weak.
The following table gives approximate values of $\delta$ for small
$a$, $b$, and $c$:
\begin{center} \begin{tabular}{|c|c|c|c|}\hline
$a$ & $b$ & $c$ & $\delta$ \\ \hline
2&3&5&0.7292\\ \hline
2&3&7&0.7407\\ \hline
2&5&7&0.8235\\ \hline
2&5&9&0.8187\\ \hline
2&7&9&0.8709\\ \hline
3&4&5&0.7093\\ \hline
3&4&7&0.7934\\ \hline
3&5&7&0.8239 \\\hline
3&5&8&0.8212 \\\hline
3&7&8&0.8727 \\\hline
\end{tabular} \end{center}
%Unfortunately, the program is very slow for higher precisions.
These results were obtained by incrementing $d$ and looking for
convergence to 4 decimal places.
We also approximated $\delta_S$ using a naive algorithm (based on
Lemma~1) for large $n$.
Numerical convergence occurred at values of $d$ slightly lower than
the bound given above. %, suggesting it could be improved.

%%%%%%%%%%%%%%%%%%%%%%%%%%%%%%%%%%%%%%%%%%%%%%%%%%%%%%%%%%%%%%
%\bibliographystyle{plainnat}
%\bibliography{abmult}
%\bibliographystyle{myNatbibStyle}
%\bibliography{myBibliography}
%%%%%%%%%%%%%%%%%%%%%%%%%%%%%%%%%%%%%%%%%%%%%%%%%%%%%%%%%%%%%%

\def\soft#1{\leavevmode\setbox0=\hbox{h}\dimen7=\ht0\advance \dimen7
  by-1ex\relax\if t#1\relax\rlap{\raise.6\dimen7
  \hbox{\kern.3ex\char'47}}#1\relax\else\if T#1\relax
  \rlap{\raise.5\dimen7\hbox{\kern1.3ex\char'47}}#1\relax \else\if
  d#1\relax\rlap{\raise.5\dimen7\hbox{\kern.9ex \char'47}}#1\relax\else\if
  D#1\relax\rlap{\raise.5\dimen7 \hbox{\kern1.4ex\char'47}}#1\relax\else\if
  l#1\relax \rlap{\raise.5\dimen7\hbox{\kern.4ex\char'47}}#1\relax \else\if
  L#1\relax\rlap{\raise.5\dimen7\hbox{\kern.7ex
  \char'47}}#1\relax\else\message{accent \string\soft \space #1 not
  defined!}#1\relax\fi\fi\fi\fi\fi\fi}

\end{document}